\documentclass[12pt]{article}
\usepackage{latexsym, amssymb}

\DeclareMathAlphabet\gothic{U}{euf}{m}{n}

\makeatletter
\def\eqnarray{\stepcounter{equation}\let\@currentlabel=\theequation
\global\@eqnswtrue
\tabskip\@centering\let\\=\@eqncr
$$\halign to \displaywidth\bgroup\hfil\global\@eqcnt\z@
  $\displaystyle\tabskip\z@{##}$&\global\@eqcnt\@ne
  \hfil$\displaystyle{{}##{}}$\hfil
  &\global\@eqcnt\tw@ $\displaystyle{##}$\hfil
  \tabskip\@centering&\llap{##}\tabskip\z@\cr}

\def\endeqnarray{\@@eqncr\egroup
      \global\advance\c@equation\m@ne$$\global\@ignoretrue}

\def\@yeqncr{\@ifnextchar [{\@xeqncr}{\@xeqncr[5pt]}}
\makeatother

\begin{document}
\bibliographystyle{tom}

\newtheorem{lemma}{Lemma}[section]
\newtheorem{thm}[lemma]{Theorem}
\newtheorem{cor}[lemma]{Corollary}
\newtheorem{voorb}[lemma]{Example}
\newtheorem{rem}[lemma]{Remark}
\newtheorem{prop}[lemma]{Proposition}
\newtheorem{stat}[lemma]{{\hspace{-5pt}}}
\newtheorem{obs}[lemma]{Observation}
\newtheorem{defin}[lemma]{Definition}

\newenvironment{remarkn}{\begin{rem} \rm}{\end{rem}}
\newenvironment{exam}{\begin{voorb} \rm}{\end{voorb}}
\newenvironment{defn}{\begin{defin} \rm}{\end{defin}}
\newenvironment{obsn}{\begin{obs} \rm}{\end{obs}}

\newenvironment{emphit}{\begin{itemize} }{\end{itemize}}

\newcommand{\gota}{\gothic{a}}
\newcommand{\gotb}{\gothic{b}}
\newcommand{\gotc}{\gothic{c}}
\newcommand{\gote}{\gothic{e}}
\newcommand{\gotf}{\gothic{f}}
\newcommand{\gotg}{\gothic{g}}
\newcommand{\gothh}{\gothic{h}}
\newcommand{\gotk}{\gothic{k}}
\newcommand{\gotm}{\gothic{m}}
\newcommand{\gotn}{\gothic{n}}
\newcommand{\gotp}{\gothic{p}}
\newcommand{\gotq}{\gothic{q}}
\newcommand{\gotr}{\gothic{r}}
\newcommand{\gots}{\gothic{s}}
\newcommand{\gotu}{\gothic{u}}
\newcommand{\gotv}{\gothic{v}}
\newcommand{\gotw}{\gothic{w}}
\newcommand{\gotz}{\gothic{z}}
\newcommand{\gotA}{\gothic{A}}
\newcommand{\gotB}{\gothic{B}}
\newcommand{\gotG}{\gothic{G}}
\newcommand{\gotL}{\gothic{L}}
\newcommand{\gotS}{\gothic{S}}
\newcommand{\gotT}{\gothic{T}}

\newcommand{\mn}{\marginpar{\hspace{1cm}*} }
\newcommand{\mnn}{\marginpar{\hspace{1cm}**} }

\newcommand{\mnq}{\marginpar{\hspace{1cm}*???} }
\newcommand{\mnnq}{\marginpar{\hspace{1cm}**???} }

\newfont{\fontcmrten}{cmr10}
\newcommand{\slbrl}{\mbox{\fontcmrten (}}
\newcommand{\slbrr}{\mbox{\fontcmrten )}}

\newcounter{teller}
\renewcommand{\theteller}{\Roman{teller}}
\newenvironment{tabel}{\begin{list}%
{\rm \bf \Roman{teller}.\hfill}{\usecounter{teller} \leftmargin=1.1cm
\labelwidth=1.1cm \labelsep=0cm \parsep=0cm}
                      }{\end{list}}

\newcounter{tellerr}
\renewcommand{\thetellerr}{(\roman{tellerr})}
\newenvironment{subtabel}{\begin{list}%
{\rm  (\roman{tellerr})\hfill}{\usecounter{tellerr} \leftmargin=1.1cm
\labelwidth=1.1cm \labelsep=0cm \parsep=0cm}
                         }{\end{list}}
\newenvironment{ssubtabel}{\begin{list}%
{\rm  (\roman{tellerr})\hfill}{\usecounter{tellerr} \leftmargin=1.1cm
\labelwidth=1.1cm \labelsep=0cm \parsep=0cm \topsep=1.5mm}
                         }{\end{list}}

\newcommand{\Ni}{{\bf N}}
\newcommand{\Ri}{{\bf R}}
\newcommand{\Ci}{{\bf C}}
\newcommand{\Ti}{{\bf T}}
\newcommand{\Zi}{{\bf Z}}
\newcommand{\Fi}{{\bf F}}

\newcommand{\proof}{\mbox{\bf Proof} \hspace{5pt}} 
\newcommand{\remark}{\mbox{\bf Remark} \hspace{5pt}}
\newcommand{\ruimte}{\vskip10.0pt plus 4.0pt minus 6.0pt}

\newcommand{\simh}{{\stackrel{{\rm cap}}{\sim}}}
\newcommand{\ad}{{\mathop{\rm ad}}}
\newcommand{\Ad}{{\mathop{\rm Ad}}}
\newcommand{\Aut}{\mathop{\rm Aut}}
\newcommand{\arccot}{\mathop{\rm arccot}}
\newcommand{\capp}{{\mathop{\rm cap}}}
\newcommand{\rcapp}{{\mathop{\rm rcap}}}
\newcommand{\diam}{\mathop{\rm diam}}
\newcommand{\divv}{\mathop{\rm div}}
\newcommand{\codim}{\mathop{\rm codim}}
\newcommand{\RRe}{\mathop{\rm Re}}
\newcommand{\IIm}{\mathop{\rm Im}}
\newcommand{\Tr}{{\mathop{\rm Tr}}}
\newcommand{\Vol}{{\mathop{\rm Vol}}}
\newcommand{\card}{{\mathop{\rm card}}}
\newcommand{\supp}{\mathop{\rm supp}}
\newcommand{\sgn}{\mathop{\rm sgn}}
\newcommand{\essinf}{\mathop{\rm ess\,inf}}
\newcommand{\esssup}{\mathop{\rm ess\,sup}}
\newcommand{\Int}{\mathop{\rm Int}}
\newcommand{\Leibniz}{\mathop{\rm Leibniz}}
\newcommand{\lcm}{\mathop{\rm lcm}}
\newcommand{\loc}{{\rm loc}}

\newcommand{\mod}{\mathop{\rm mod}}
\newcommand{\spann}{\mathop{\rm span}}
\newcommand{\one}{1\hspace{-4.5pt}1}

\newcommand{\DWR}{}

\hyphenation{groups}
\hyphenation{unitary}

\newcommand{\tfrac}[2]{{\textstyle \frac{#1}{#2}}}

\newcommand{\cb}{{\cal B}}
\newcommand{\cc}{{\cal C}}
\newcommand{\cd}{{\cal D}}
\newcommand{\ce}{{\cal E}}
\newcommand{\cf}{{\cal F}}
\newcommand{\ch}{{\cal H}}
\newcommand{\ci}{{\cal I}}
\newcommand{\ck}{{\cal K}}
\newcommand{\cl}{{\cal L}}
\newcommand{\cm}{{\cal M}}
\newcommand{\cn}{{\cal N}}
\newcommand{\co}{{\cal O}}
\newcommand{\cs}{{\cal S}}
\newcommand{\ct}{{\cal T}}
\newcommand{\cx}{{\cal X}}
\newcommand{\cy}{{\cal Y}}
\newcommand{\cz}{{\cal Z}}

\newcommand{\wtozp}{W^{1,2}\raisebox{10pt}[0pt][0pt]{\makebox[0pt]{\hspace{-34pt}$\scriptstyle\circ$}}}
\newlength{\hightcharacter}
\newlength{\widthcharacter}
\newcommand{\covsup}[1]{\settowidth{\widthcharacter}{$#1$}\addtolength{\widthcharacter}{-0.15em}\settoheight{\hightcharacter}{$#1$}\addtolength{\hightcharacter}{0.1ex}#1\raisebox{\hightcharacter}[0pt][0pt]{\makebox[0pt]{\hspace{-\widthcharacter}$\scriptstyle\circ$}}}
\newcommand{\cov}[1]{\settowidth{\widthcharacter}{$#1$}\addtolength{\widthcharacter}{-0.15em}\settoheight{\hightcharacter}{$#1$}\addtolength{\hightcharacter}{0.1ex}#1\raisebox{\hightcharacter}{\makebox[0pt]{\hspace{-\widthcharacter}$\scriptstyle\circ$}}}
\newcommand{\scov}[1]{\settowidth{\widthcharacter}{$#1$}\addtolength{\widthcharacter}{-0.15em}\settoheight{\hightcharacter}{$#1$}\addtolength{\hightcharacter}{0.1ex}#1\raisebox{0.7\hightcharacter}{\makebox[0pt]{\hspace{-\widthcharacter}$\scriptstyle\circ$}}}

 \thispagestyle{empty}

\begin{center}
{\Large{\bf  Gru\v{s}in operators, Riesz transforms  
 }}\\[4mm] 
{\Large{\bf  and  nilpotent Lie groups}}  \\[4mm]
\large  Derek W. Robinson$^1$ and Adam Sikora$^2$\\[2mm]

\normalsize{February 2014}
\end{center}

\vspace{5mm}

\begin{center}
{\bf Abstract}
\end{center}

\begin{list}{}{\leftmargin=1.8cm \rightmargin=1.8cm \listparindent=10mm 
   \parsep=0pt}
   \item

We establish that the Riesz transforms of all orders corresponding to the Gru\v{s}in operator $H_N=-\nabla_{x}^2-|x|^{2N}\,\nabla_{y}^2$,
 and the first-order operators $(\nabla_{x},x^\nu\,\nabla_{y})$
where $x\in \Ri^n$, $y\in\Ri^m$,  $N\in\Ni_+$, and $\nu\in\{1,\ldots,n\}^N$, are bounded on $L_p(\Ri^{n+m})$ for all $p\in\langle1,\infty\rangle$ and are also weak-type $(1,1)$.
Moreover, the transforms of order less than or equal to $N+1$ corresponding to $H_N$ and 
the operators  $(\nabla_{x}, |x|^N\nabla_{y})$ are bounded on $L_p(\Ri^{n+m})$ for all $p\in\langle1,\infty\rangle$. 
But all  transforms of order $N+2$  are bounded if  and only if $p\in\langle1,n\rangle$.
The proofs are based on the observation that the $(\nabla_{x},x^\nu\,\nabla_{y})$ generate a finite-dimensional nilpotent Lie algebra,
the corresponding connected, simply connected, nilpotent Lie group is isometrically represented on the spaces $L_p(\Ri^{n+m})$
and  $H_N$ is the corresponding sublaplacian

\vfill

\end{list}

\vspace{1cm}


\vspace{1cm}

\vspace{1cm}

\noindent
\begin{tabular}{@{}cl@{\hspace{10mm}}cl}
1. & Centre for Mathematics & 
  2. & Department of Mathematics  \\
&\hspace{15mm} and its Applications  & 
  &Macquarie University  \\
& Mathematical Sciences Institute & 
  & Sydney, NSW 2109  \\
& Australian National University& 
  & Australia \\
& Canberra, ACT 0200  & {}
  & \\
& Australia & {}
  & \\
& derek.robinson@anu.edu.au & {}
  &sikora@ics.mq.edu.au \\
\end{tabular}

\newpage
\setcounter{page}{1}

\section{Introduction}\label{S1}

The primary aim of this note is to prove the boundedness of the Riesz transforms  corresponding to the Gru\v{s}in
operators 
\begin{equation}
H_N=-\nabla_{x}^2-|x|^{2N}\,\nabla_{y}^2\;,
\label{egrrt1.1}
\end{equation}
and the family of first-order operators
$(\nabla_x,\,x^\nu\,\nabla_y)$ where $x\in\Ri^n$, $y\in\Ri^m$, $N\in\Ni_+$ and  the multi-index $\nu$ takes all values in the set $\{1,\ldots,n\}^N$.
We establish that  the transforms of arbirtrary order are bounded on the spaces $L_p(\Ri^{n+m})$ for all $p\in\langle1,\infty\rangle$ and are also weak-type $(1,1)$.
The proof of boundedness is based on the observation that there is an underlying symmetry which  arises because the operators
$(\nabla_x,\,x^\nu\,\nabla_y)$ generate a finite-dimensional nilpotent Lie algebra.
This allows one to identify the transforms with the Riesz transforms associated with an isometric representation of a nilpotent Lie group and an affiliated subelliptic operator.
It is known, however, that the Riesz transforms corresponding to very general subcoercive operators in the left, or right, regular representation of a nilpotent Lie group are bounded \cite{BER} \cite{ERS}.
We will show that one can adapt the Coifman-Weis transference techniques \cite{CW} to deduce from this result that the transforms in a every isometric representation of the group  are also  bounded.
As an immediate corollary one deduces that the Riesz transforms connected to the Gru\v{s}in operators are bounded.
This result was established recently by Jotsaroop, Sanjay and Thangavelu \cite{JST} for the special case $m=1=N$ using quite different arguments (see also \cite{ABCS}).

The symmetry techniques that we exploit are sensitive to the choice of the first-order operators
$(\nabla_x,\,x^\nu\,\nabla_y)$ and are also dependent on the assumption that $N$ is an integer. 
In order to gain some insight into the  difficulties that arise in the non-integer situation we next examine an alternative
choice of transforms.
We consider the transforms corresponding to $H_N$ and the first-order operators $(\nabla_x,|x|^N\nabla_y)$.
These operators no longer generate a Lie algebra.
Nevertheless, if $n\geq 2$  one can   deduce from the group-theoretic results  that the corresponding Riesz transforms of all orders less than or equal to $N+1$  are bounded on $L_p(\Ri^{n+m})$ with $p\in\langle1,\infty\rangle$
and also deduce that these transforms are weak-type $(1,1)$. 
In addition one can demonstrate that some  of the transforms of order $N+2$ are unbounded on  the $L_p$-spaces 
with $p\geq n$
but all the transforms of order $N+2$ are bounded on the $L_p$-spaces with $p\in\langle1,n\rangle$.

In Section~\ref{S2} we describe briefly the necessary background on Riesz transforms on nilpotent Lie groups and 
establish a very general boundedness result for isometric group representations.
Theorem~\ref{tgrrt2.1} and its corollaries provide the basis for the discussion of the Gru\v{s}in  transforms
in Section~\ref{S3} although the details of their proofs are largely superfluous for comprehension of the Gru\v{s}in analysis.

\section{Riesz transforms and nilpotent groups}\label{S2}

Let $G$ be a $d$-dimensional connected nilpotent Lie group with
Haar measure $dg$ and  Lie algebra $\gotg$.
Further let  $r$ denote   the rank of $\gotg$. 
Fix an algebraic basis  $a_1,\ldots,a_{d'}$  and 
a corresponding vector space basis $a_1,\ldots,a_{d'}, a_{d'+1},\ldots, a_d$  of $\gotg$.
The group acts isometrically by left  translations $L$
on the spaces $L_p(G)=L_p(G\,;dg)$ for all $p\in[1,\infty]$. 
Next let $A_i=dL(a_i)$  denote the corresponding generators of 
left translations, i.e.\ $A_i$ is the generator of the one-parameter group
 $u_i\in \Ri \mapsto L(\exp(u_ia_i))$, where $u\in\Ri^d\mapsto \exp(u.a)\in G$ denotes the usual exponential map.
 We use multi-index notation for products of the generators corresponding to the algebraic basis.
For example,  if $\alpha=(i_1,\ldots, i_n)$ with $i_j\in\{1,\ldots, d'\}$ then $A^\alpha=A_{i_1}\ldots A_{i_n}$ and $|\alpha|=n$.
We also set $L^{\,\prime}_{p;n}(G)=\bigcap_{\{\alpha:|\alpha|\leq n\}}D(A^\alpha)$ with the usual graph norm.

Since we wish to apply various results from \cite{ERS} it is also convenient to introduce the nilpotent Lie  group 
$\widetilde G$ with $d'$ generators which is free of step $r$.
First define $\tilde\gotg $ as the quotient of the free Lie algebra with $d'$ generators 
${\tilde a}_1,\ldots,{\tilde a}_{d'}$ by the ideal generated by the commutators of order at least $r+1$. 
Then  $\widetilde G$ is the connected simply connected 
Lie group with Lie algebra $\tilde\gotg $.
It is non-compact,  stratified and acts
 isometrically by left  translations $\tilde L$ on the spaces $L_p(\widetilde G)=L_p(\widetilde G\,;d\tilde g)$
where $d\tilde g$ denotes the Haar measure of $\widetilde G$.

Next let $X=(X,\mu)$ be a $\sigma$-finite measure space and $U$  an  isometric 
representation of $G$ on the spaces $L_p(X)=L_p(X\,;\mu)$.
Further let  $D_k=dU(a_k)$   denote the  corresponding group generators,
$D^\beta$ their products  and set  $L^{\,\prime}_{p;n}(X)=\bigcap_{\{\beta:|\beta|\leq n\}}D(D^\beta)$.
Then  for  each $k\in L_1(G)$ we define bounded  operators $K_L$ and $K_U$ acting on $L_p(G)$ and $L_p(X)$, respectively by
\[
K_L=\int_Gdg\, k(g)\,L(g)\;\;\;\;\;{\rm and}\;\;\;\;\;K_U=\int_Gdg\, k(g)\,U(g)
\;.
\]
Note that the $K_L$ act by convolution, $K_L\varphi=k*\varphi$, where $(k*\varphi)(g)=\int_Gdh\,k(g)\varphi(h^{-1}g)$.

\smallskip

One can transfer estimates on  $K_L$ to $K_U$   by the Coifman--Weiss transference theorem.

\begin{prop}\label{pgrrt2.1} If $k\in L_1(G)$ then
$\|K_U\|_{L_p(X)\to L_p(X)}\leq \|K_L\|_{L_p(G)\to L_p(G)}$
for all $p\in [1,\infty]$.
\end{prop}

If $k$ has compact support then the result is given by Theorem~2.4 of \cite{CW}.
The general  result follows by approximating $k$  with compactly supported functions
and using the isometry of  the group actions $U$ and $L$.\smallskip

Next observe that \vspace{-2mm}
\begin{eqnarray*}
U(\exp(u.a))K_U&=&\int_Gdg\, k(g)\,U(\exp(u.a)g)\\[5pt]
&=&\int_Gdg\, k(\exp(-u.a)g)\,U(g)=\int_Gdg\, (L(\exp(u.a))k)(g)\,U(g)
\end{eqnarray*}
by  the group property and invariance of Haar measure. 
Therefore if $k\in L^{\,\prime}_{1;n}(G)$ then $K_UL_p(X)\subseteq L'_{p;n}(X)$ and
\begin{equation}
D^\alpha K_U=\int_G dg\,(A^\alpha k)(g)U(g)
\label{egrrt2.4}
\end{equation}
for  all $\alpha $ with $|\alpha|\leq n$.

Let $p\in\langle1,\infty\rangle$ and let $m$ be an even positive integer.
Following \cite{ERS} we introduce  homogeneous $m$-th order operators
\[
H = \sum_{|\alpha| = m} c_\alpha \, A^\alpha
\;\;\;\;\;\;\;{\rm and}  \;\;\;\;\;\;\;\widetilde H = \sum_{|\alpha| = m} c_\alpha \, \tilde A^\alpha
\;,
\]
with $c_\alpha \in\Ci$, on $L_{p;m}^{\,\prime}(G)$ and $L_{p;m}^{\,\prime}(\widetilde G)$, respectively.
These operators are closable and for simplicity we assume that $H$ and $\widetilde H$ denote the closures.
Further we assume  $H$ is subcoercive of step $r$ in the sense  of \cite{ER7}, i.e.\ 
$\widetilde H$ satisfies the G\aa rding inequality 
\[
\RRe (\tilde \varphi, \widetilde H \tilde \varphi)
\geq \mu \sum_{|\alpha| = m/2} \|\tilde A^\alpha \tilde \varphi\|_{\tilde 2}
\]
for some $\mu>0$ uniformly for all $\tilde \varphi \in C_c^\infty(\widetilde G)$.

\begin{remarkn}\label{rgrrt2.1}
The subcoercivity requirement for second-order operators with $r \geq 2$  is equivalent to the assumption that 
$\Re C = 2^{-1} (C + C^*)\geq \mu I>0$, where $C$ is the 
 $d' \times d'$-matrix given by
$C_{ij} = - c_{(i,j)}$ (see \cite{ER8}, Proposition 3.7).
In particular the sublaplacian $-\sum^{d'}_{j=1}A_j^2$ is subcoercive.
\end{remarkn}

It follows  that each subcoercive operator
$H$  generates a holomorphic semigroup $S_t$
 with a $C^\infty$-kernel $k_t$ on $L_p(G)$.
Specifically 
\[
S_t=\int_G dg\,k_t(g)\,L(g)
\]
for all $t>0$ or, in the previous notation, $S_t=k_{t, L}$.
(We refer to  \cite{ER7} and \cite{ER8} for details on the semigroup and its kernel.)
Next  for each $t>0$ we define  the operator 
\[
S^U_t=\int_G dg\,k_t(g)\,U(g)\;\;(\;=k_{t, U}\,)
\]
on $L_p(X)$.
Since $k_t$ is a convolution semigroup and $U$ a group representation  it follows that  $S^U$  is a  
semigroup.
Moreover, $\|S^U_t\|_{L_p(X)\to L_p(X)}\leq \|S_t\|_{L_p(G)\to L_p(G)}$ by Proposition~\ref{pgrrt2.1}.
Now let
\[
H_D=\sum_{|\alpha| = m} c_\alpha \, D^\alpha
\]
on $L^{\,\prime}_{p;m}(X)$.
It follows from (\ref{egrrt2.4}) that 
$S^U_tL_p(X)\subseteq L^{\,\prime}_{p;m}(X)$ for all $t>0$  and 
\[
H_DS^U_t=\int_Gdg\,(Hk_t)(g)\,U(g)
=-{{d}\over{dt}}S^U_t
\]
for  all $t>0$.
Therefore the closure of $H_D$, which for simplicity we also denote by $H_D$,  is the generator of $S^U$.
Moreover, $\|H_DS^U_t\|_{L_p(X)\to L_p(X)}\leq \|HS_t\|_{L_p(G)\to L_p(G)}$ by another application of Proposition~\ref{pgrrt2.1}.
Hence $S^U$ inherits the holomorphy properties of $S$.

Now we are prepared to prove boundedness of the Riesz transforms corresponding to the derivatives $D$ and the operator $H_D$.
These transforms are formally given  by  $R^U_\alpha=D^\alpha H_D^{-|\alpha|/m}$ on $L_p(X)$
but it is not clear that the operators are even densely-defined. 
In fact the operators are given by integral  kernels which  are logarithmically divergent both at the identity and at infinity.
These difficulties can, however, be overcome by the techniques used in earlier papers \cite{BER} \cite{ERS} to handle the Riesz  transforms $R_\alpha$ and $\widetilde R_\alpha$ associated with $H$ and $\widetilde H$ on $L_p(G)$  and $L_p(\widetilde G)$, respectively.
Theorem~4.4 of \cite{ERS}  establishes that the $R_\alpha$ and $\widetilde R_\alpha$ extend to bounded operators
but the same is true for the~$R^U_\alpha$.

\begin{thm}\label{tgrrt2.1}
The Riesz transforms $R^U_\alpha$ extend to bounded operators on $L_p(X)$ for each
$p\in\langle1,\infty\rangle$.
Moreover, there exist $a_p,b_p>0$ 
such that 
\begin{equation}\|R^U_\alpha\|_{L_p(X)\to L_p(X)}
\leq a_p\,b_p^{\,k}\,\| R_\alpha\|_{L_p( G)\to L_p( G)}
\label{egrrt2.50}
\end{equation}
where 
 $k$ is the integer part of $|\alpha|/m$.
\end{thm}
\proof\
The theorem is a direct analogue of Theorem~4.4 of \cite{ERS}  and it follows by 
 an adaptation of the arguments of \cite{BER} \cite{ERS} supplemented by the transference 
statement of Proposition~\ref{pgrrt2.1}.
But first note that if $G$ is compact then we define the transforms $R_\alpha$ to be zero on the constant functions and analyze the transforms on $L_p(G)$ modulo the constant functions.

The proof of the theorem is based on estimation of the regularized transforms
\[
R_{\alpha;\nu,\varepsilon}=A^\alpha (\nu I+H)^{-|\alpha|/m}(I+\varepsilon H)^{-n}
\]
on $L_p(G)$ with $\varepsilon, \nu>0$ and $n$ a large positive integer together with the corresponding regularizations
$R^U_{\alpha;\nu,\varepsilon}$ of $R^U_\alpha$  on 
$L_p(X)$.

The first point  is that  $R_{\alpha;\nu,\varepsilon}$ and $R^U_{\alpha;\nu,\varepsilon}$  are given by a kernel 
$k_{\alpha;\nu,\varepsilon}\in L_1(G)$.
The introduction  of the $\nu$-term ensures integrability at infinity and the factor $(I+\varepsilon H)^{-n}$ ensures integrability at the identity.
Consequently 
\begin{equation}
\|R^U_{\alpha;\nu,\varepsilon}\|_{L_p(X)\to L_p(X)}\leq \|R_{\alpha;\nu,\varepsilon}\|_{L_p(G)\to L_p(G)}
\label{egrrt2.5}
\end{equation}
by Proposition~\ref{pgrrt2.1}.
The second point is that since the $R_\alpha$ extend to bounded operators on $L_p(G)$ one can estimate the right hand norm
uniformly for $\nu,\varepsilon\in \langle0,1]$ by
 a variation of the argument of \cite{ERS}.

The starting point is the observation that 
\[
\|R_{\alpha;\nu,\varepsilon}\|_{p\to p}\leq \|R_{\alpha}\|_{p\to p}\,\|H^{|\alpha|/m}(\nu I+H)^{-|\alpha|/m}
(I+\varepsilon H)^{-n}\|_{p\to p}
\]
where for  brevity we have set $\|\,\cdot\,\|_{L_p(G)\to L_p(G)}=\|\,\cdot\,\|_{p\to p}$.
But   $|\alpha|/m=k+\gamma$ with $k\in\Ni_+$ and $\gamma\in[0,1\rangle$.
Therefore
\begin{eqnarray*}
\|R_{\alpha;\nu,\varepsilon}\|_{p\to p}\leq \|R_{\alpha}\|_{p\to p}\left(\|H(\nu I+H)^{-1}\|_{p\to p}\right)^k
\|H^\gamma(\nu I+H)^{-\gamma}(I+\varepsilon H)^{-n}\|_{p\to p}
\;.
\end{eqnarray*}

Now consider the case $\gamma=0$.
First one has \vspace{-2mm}
\[
\|H(\nu I+H)^{-1}\|_{p\to p}\leq \Big\|\,\nu \int^\infty_0dt\,e^{-\nu t}(I-S_t)\Big\|_{p\to p}\leq b_p
\]
with $b_p=\sup_{t>0}\|I-S_t\|_{p\to p}$, which is  finite  because $S$ is uniformly bounded by (3) of \cite{ERS}.
Secondly,\vspace{-2mm}
\[
\|(I+\varepsilon H)^{-n}\|_{p\to p}
\leq \Gamma(n)^{-1}\int^\infty_0dt\,e^{-t}t^{-1+n}\|S_{\varepsilon t}\|_{p\to p}\leq a_p
\]
where $a_p=\sup_{t>0}\|S_t\|_{p\to p}$ which is again finite.
Combining these estimates  one has
\begin{equation}
\|R_{\alpha;\nu,\varepsilon}\|_{p\to p}\leq a_p\,b_p^{\,k}\,\|R_{\alpha}\|_{p\to p}
\label{egrrt2.51}
\end{equation}
uniformly for $\nu,\varepsilon>0$.

Next consider the case $\gamma\in\langle0,1\rangle$.
Then $H(\nu I+H)^{-1}$ is bounded and generates a uniformly bounded semigroup.
Therefore  its fractional power is given by the algorithm
\null\vspace{-3mm}
\begin{eqnarray*}
H^\gamma(\nu I+H)^{-\gamma}&=& 
n_\gamma^{-1}\int^\infty_0d\lambda\,\lambda^{-1+\gamma}(1+\lambda)^{-1}
\, H(\lambda(1+\lambda)^{-1} \nu I+ H)^{-1}
\end{eqnarray*}
where $n_\gamma=\int^\infty_0d\lambda\,\lambda^{-1+\gamma}(1+\lambda)^{-1}$.
Hence
\begin{eqnarray*}
\|H^\gamma(\nu I+H)^{-\gamma}(I+\varepsilon H)^{-n}\|_{p\to p}
&\leq &\sup_{\nu,\varepsilon>0}\|H_\nu( I+H_\nu)^{-1}(I+\varepsilon H)^{-n}\|_{p\to p}
\end{eqnarray*}
where $H_\nu=\nu^{-1}H$.
But $H_\nu$ generates the uniformly bounded semigroup  $S^\nu_t=S_{\nu^{-1}t}$.
Therefore
\null\vspace{-2mm}
\[
H_\nu(I+ H_\nu)^{-1}
(I+\varepsilon H)^{-n}=\Gamma(n)^{-1}
\int^\infty_0ds\int^\infty_0dt\,e^{-(s+t)}s^{-1+n}
S_{\varepsilon s}(I-S_{\nu^{-1}t})
\;.
\]
Consequently one deduces that 
\[
\|H_\nu( I+H_\nu)^{-1}(I+\varepsilon H)^{-n}\|_{p\to p}\leq a_p
\]
where $a_p$ is now given by $a_p=\sup_{s,t>0}\|S_s-S_t\|_{p\to p}$.
It then follows from combination of these estimates that (\ref{egrrt2.51}) is again valid with the
new choice of $a_p$.

Finally, if $\varphi\in L^{\,\prime}_{p;\infty}(X)$ then 
\begin{eqnarray*}
\|D^\alpha\varphi\|_{L_p(X)}&=&
\|R_{\alpha;\nu,\varepsilon}^U(I+\varepsilon H_D)^n(\nu I+H_D)^{|\alpha|/m}\varphi\|_{L_p(X)}\\[5pt]
&\leq&\|R_{\alpha;\nu,\varepsilon}\|_{p\to p}\,\|(I+\varepsilon H_D)^n(\nu I+H_D)^{|\alpha|/m}\varphi\|_{L_p(X)}\\[5pt]
&\leq&a_p\,b_p^k\,\|R_{\alpha}\|_{p\to p}\,\|(I+\varepsilon H_D)^n(\nu I+H_D)^{|\alpha|/m}\varphi\|_{L_p(X)}
\end{eqnarray*}
where we have used (\ref{egrrt2.5}) and (\ref{egrrt2.51}).
Therefore taking the limits $\nu,\varepsilon\to0$ one deduces that 
\[
\|D^\alpha\varphi\|_{L_p(X)}\leq a_p\,b_p^k\,\|R_{\alpha}\|_{p\to p}\,\|H_D^{|\alpha|/m}\varphi\|_{L_p(X)}
\;.
\]
But these bounds then extend to all $\varphi\in D(H_D^{|\alpha|/m})$ by continuity (see \cite{ERS}, proof of Lemma~4.2).
One immediately deduces that the Riesz transforms extend to bounded operators and that (\ref{egrrt2.50}) is valid.
\hfill$\Box$

\bigskip

The theorem has a number of corollaries similar to those given in \cite{ERS} for the operator $H$ on the spaces $L_p(G)$.

\begin{cor}\label{cgrrt2.1}
For each $n\in\Ni$ and $p\in\langle1,\infty\rangle$ one has $D(H_D^{n/m})=L^{\,\prime}_{p;n}(X)$.
Moreover, there are $c_{p,n},c^{\,\prime}_{p,n}>0$ such that
\[
c_{p,n}\,\max_{|\alpha|=n}\|D^\alpha\varphi\|_{L_p(X)}\leq \|H_D^{n/m}\varphi\|_{L_p(X)}
\leq 
c{\,\prime}_{p,n}\,\max_{|\alpha|=n}\|D^\alpha\varphi\|_{L_p(X)}
\]
for all $\varphi\in L^{\,\prime}_{p;n}(X)$.
\end{cor}
\proof\
The left hand bound is just a restatement of the boundedness of the Riesz transforms and the proof of the right 
hand bound is identical to the proof of the analogous result in Corollary~4.3 of \cite{ERS}. \hfill$\Box$

\begin{cor}\label{cgrrt2.2}
The operators $D^\alpha H_D^{-(|\alpha|+|\beta|)/m}D^\beta$ extend to bounded operators on $L_p(X)$
for each $p\in\langle1,\infty\rangle$.
\end{cor}
\proof\
This follows by the argument used to prove Corollary~4.6 of \cite{ERS}. \hfill$\Box$

\bigskip

The last  corollary could also be deduced  by transference from the result for $H$.
The key point is the identity 
\begin{equation}
D^\alpha K_U D^\beta =(-1)^{|\beta|}\int_G dg\,(A^\alpha B^{\beta_*}k)(g)U(g)
\label{egrrt2.41}
\end{equation}
where the $B_i$ are the generators of right translations of $G$ and $\beta_*$
is the multi-index obtained from $\beta$ by reversing its order.
This identity is a generalization of (\ref{egrrt2.4}) and is proved by a similar argument.
Note that the semigroup $S^U$ acting on the $L_p(X)$-spaces can be represented by a distributional kernel
kernel $(x,y)\in X\times X\mapsto K_t(x\,;y)$  and then (\ref{egrrt2.41}) gives
\[
\int_Xdy\,(D_x^\alpha D_y^\beta K_t)(x\,;y)\varphi(y)=(-1)^{|\beta|}\int_G dg\,(A^\alpha B^{\beta_*}k)(g)(U(g)\varphi)(x)
\]
where $D_x$ and $ D_y$ indicate that the operators $D$ act on the $x$ and $y$ variables respectively.

\smallskip

Finally one has a weak-type $(1,1)$ statement.

\begin{cor}\label{cgrrt2.3}
The Riesz transforms $R^U_\alpha$ extend to operators of weak-type $(1,1)$.
\end{cor}

This also follows by transference, although it can be deduced in various other ways. 
Proposition~4.7 of \cite{ERS} establishes by standard singular integration arguments (see \cite{Ste3}, Chapter~1) that the $R_\alpha$ 
extend to operators of weak-type $(1,1)$.
Then the corollary follows from the second transference theorem of Coifman and Weiss, \cite{CW} Theorem~2.6.
We omit the details.

\section{Gru\v{s}in operators}\label{S3}

In this section we apply
Theorem~\ref{tgrrt2.1}  to the Gru\v{s}in operator defined by (\ref{egrrt1.1})  in the introduction.
But we begin by filling in some of the preliminary details outlined in Section~\ref{S1}.

If  $\nu=(\nu_1, \ldots, \nu_n)\in \Ni_+^n$ is an $n$-tuple of positive integers 
we set  $x^\nu=x_1^{\nu_1}\ldots x_{n}^{\nu_{n}}$ for each $x\in\Ri^n$.
Further let $I_n(N)=\{\nu:|\nu|=N\}$ where  $|\nu|=\nu_1+\ldots+\nu_n=N$ and set $\nu!=\nu_1!\ldots\nu_n!$.
The binomial expansion gives
\[
|x|^{2N}
=\sum_{\nu\in I_n(N)} (a_\nu\, x^\nu)^2
\]
with $a_\nu=(N!/\nu!)^{1/2}$.
Therefore
\[
H_N=-\sum^n_{j=1}X_j^2-\sum^m_{k=1} \sum_{\nu\in I_n(N)} X_{\nu,k}^2
\]
with  $X_j=\partial_{x_j}$ and $X_{\nu,k}=a_\nu\,(x^\nu\,\partial_{y_k})$.
Thus   $H_N$ is the `sum of squares' of a finite family of first-order operators
$\cx_n\bigcup \cx_{N,m}$ where 
\[
\cx_n=\{X_j:1\leq j\leq n\}\;\;\;\;\;{\rm and}\;\;\;\;\;\cx_{N,m}=\{X_{\nu,k}:\nu\in I_n(N)\,,\, 1\leq k\leq m\}
\;.
\]
Both sets are abelian but
$
[\,X_j, Y_{\nu,k}\,]=a_\nu\,(\partial_{x_j}x^\nu)\,\partial_{y_k}\in \cy_{N-1,m}.
$
Therefore $\cx_n\bigcup \cx_{N,m}$ generates a nilpotent Lie algebra $\gotg_N$ of rank $N+1$.
(The generating set has dimension 
$d'=n+m{{N+n-1}\choose{n-1}}$
and the dimension of the Lie algebra is given by $d=n+m{{N+n}\choose{n}}$.)
It is now notationally convenient to fix an enumeration of the set $\cx_{N,m}$ and relabel it as $X_j$ with $n+1\leq j\leq d'$.
Then \vspace{-2mm}
\[
H_N=-\sum^{d'}_{j=1}X_j^2
\]
and the set $X_1,\ldots, X_{d'}$ is an algebraic basis of the nilpotent Lie algebra $\gotg_N$.

Let $G_N$ denote the connected simply connected nilpotent Lie group which has $\gotg_N$ as Lie algebra.
Since the closure of each  $X_j$  generates a one-parameter group of translations on the spaces $L_p(\Ri^{n+m})$ 
 the  family of $X_j$ generate an isometric representation of $G$ on $L_p(\Ri^{n+m})$ for each $p\in[1,\infty]$.
The Gru\v{s}in operator $H_N$  is the  corresponding sublaplacian.
In particular $H_N$ is a second-order subcoercive operator (see Remark~\ref{rgrrt2.1}).
Next 
if $s\in\Ni_+$ and $\alpha=(i_1,\ldots,i_s)\in\{1, \ldots, d'\}^s$ we set $X^\alpha=X_{i_1}\ldots X_{i_s}$ and $|\alpha|=s$.
 Then one immediately deduces the following.

\begin{thm}\label{tgrrt3.1} The Riesz transforms $X^\alpha H_N^{-|\alpha|/2}$ extend to bounded operators on each of 
the spaces $L_p(\Ri^{n+m})$ with $p\in\langle1,\infty\rangle$ and are of weak-type $(1,1)$.
\end{thm}

This is a direct corollary of Theorem~\ref{tgrrt2.1}, Corollary~\ref{cgrrt2.3} and the foregoing observations.

\medskip

Next consider the Riesz transforms associated with the
Gru\v{s}in operator $H_N $ and the $n+m$ operators~$Y_j$
 given by 
$Y_j=\partial_{x_j}$ if $j\in\{1,\ldots,n\}$ and $Y_j=|x|^N\,\partial_{y_l}$ if $j=n+l$ with $l\in\{1,\ldots,m\}$.
Now  \vspace{-2mm}
\[
H_N=-\sum^{n+m}_{j=1}Y_j^2
\]
but the $Y_j$  no longer generate a finite-dimensional Lie algebra.
For example, $[\,Y_j, Y_{n+l}\,]=(x_j/|x|)\,N\,|x|^{N-1}\partial_{y_l}$.
Nevertheless one has the following partial conclusion about the corresponding Riesz transforms.
(In the sequel the  multi-index $\alpha$ corresponds to indices in the set $\{1,\ldots, n+m\}$.)

\begin{thm}\label{tgrrt3.2} If $n\geq 2$ the Riesz transforms $Y^\alpha H_N^{-|\alpha|/2}$ with $|\alpha|\leq N+1$ extend to bounded operators on $L_p(\Ri^{n+m})$
for  each $p\in\langle1,\infty\rangle$ and are of weak-type $(1,1)$.
\end{thm}
\proof\ First consider the case $|\alpha|=1$.
If $j\leq n$ then $Y_j=X_j$.
But if  $j>n$ then 
\[
\|Y_j\varphi\|_p=\|\, |x|^N\,\partial_{y_j}\varphi\|_p\leq \|\, |x|_1^N\,\partial_{y_j}\varphi\|_p
\]
where $|x|_1$ is the $l_1$-norm of $x$.
Then, however, one has an estimate
\[
\|Y_j\varphi\|_p\leq c_N\sup_{\nu=I_n(N)}\| \,x^\nu\,\partial_{y_j}\varphi\|_p
\leq c_N \sup_{k>n}\|X_k\varphi\|_p
\;.\]

\vspace{-2mm}

\noindent
It then  follows from  Theorem~\ref{tgrrt3.1}
and Corollary~\ref{cgrrt2.1} that $D(H_N^{1/2})=\bigcap^{n+m}_{j=1}D(Y_j)$
and $\|Y_j\varphi\|_p\leq c_N\sup_{1\leq k\leq d'}\|X_k\varphi\|_p\leq c_N'\|H_N^{1/2}\varphi\|_p$
for all $j\in\{1,\dots, n+m\}$, $\varphi\in D(H^{1/2}_N)$ and $p\in\langle1,\infty\rangle$.
Thus the first-order transforms $Y_jH_N^{-1/2}$ extend to bounded operators on the $L_p$-spaces 
and are weak-type $(1,1)$ by Theorem~\ref{tgrrt3.1} and Corollary~\ref{cgrrt2.3}.

Next we proceed by induction.
We make the induction hypothesis that 
\begin{equation}
\sup_{\{\alpha:|\alpha|=M\}}
\|Y^\alpha\varphi\|_p\leq c_N\sup_{\{\beta:|\beta|=M\}}\|X^\beta\varphi\|_p
\label{egrrt3.1}
\end{equation}
for all $\varphi\in C_c^\infty(\Ri^{n+m})$ and some $M\leq N$.
It follows from the foregoing that the hypothesis is valid for $M=1$ and we use this to argue that it is also valid for
$M+1$.

Let $Y^\alpha=Y^\gamma Y_j$ with $|\alpha|=M+1$.
If $j<n$ then $Y_j=X_j$ and 
\[
\|Y^\alpha\varphi\|_p=\|Y^\gamma Y_j \varphi\|_p
\leq c_N\sup_{\{\gamma:|\gamma|=M\}}\|X^\gamma X_j\varphi\|_p
\leq 
c_N\sup_{\{\beta:|\beta|=M+1\}}\|X^\beta\varphi\|_p
\]
and (\ref{egrrt3.1}) is satisfied.
Therefore assume that $n<j\leq n+m$.
Then
\begin{eqnarray}
\|Y^\alpha\varphi\|_p
&\leq &c_N\sup_{\{\gamma:|\gamma|=M\}}\|X^\gamma Y_j\varphi\|_p\nonumber\\[5pt]
&\leq&c_N\sup_{\{\gamma:|\gamma|=M\}}\|Y_jX^\gamma \varphi\|_p
+c_N\sup_{\{\gamma:|\gamma|=M\}}\|(\ad\, X^\gamma)( Y_j)\varphi\|_p
\;.\label{egrrt3.10}
\end{eqnarray}
The first-order estimate then gives
 an upper bound $c^{\,\prime}_N\sup_{\{\beta:|\beta|=M+1\}}\|X^\beta\varphi\|_p$ on the first term on the right
as required.
It remains to estimate the second term.

The starting point is the multi-commutator algorithm
\begin{equation}
(\ad\, X^\gamma)( Y_j)=\sum_{\delta\subseteq \gamma}
(\ad\, X)^{\delta}( Y_j)X^{\gamma\backslash\delta}
\;.
\label{egrrt3.2}
\end{equation}
\vspace{-4mm}

\noindent The sum is over all non-empty ordered subsets $\delta$ of $\gamma$ and $\gamma\backslash \delta$ is the ordered set obtained from $\gamma$ by omitting the elements of $\delta$.
Now by assumption $Y_j=|x|^N\partial_{y_j}$ and each $X_j$ is either of the form $x^\nu \partial_{y_k}$ or $\partial_{x_l}$.
But each $x^\nu \partial_{y_k}$ commutes with $Y_j$ and with each $(\ad\, \partial_x)^{\delta'}(Y_j)$.
The latter observation follows because $(\ad\, \partial_x)^{\delta'}(Y_j)=( \partial_x^{\,\delta'}|x|^N)\,\partial_{y_j}$.
Therefore the only non-zero terms on the right hand side of (\ref{egrrt3.2})  correspond to multi-indices  $\delta'\subseteq \gamma$ with entries from the set $\{1,\ldots, n\}$.
Therefore
\[
(\ad\, X^\gamma)( Y_j)=\sum_{\delta'\subseteq \gamma}
(\ad\,X)^{\delta'}(Y_j)\,X^{\gamma\backslash\delta'}=
\sum_{\delta'\subseteq \gamma}
(\partial_x^{\,\delta'}|x|^N)\,\partial_{y_j}\,X^{\gamma\backslash\delta'}
\;.
\]
Now since $|\delta'|\leq |\gamma|\leq N$ and $n\geq2$  there are $a_N,b_N>0$ such that
\[
|(\partial_x^{\,\delta'}|x|^N)|\leq a_N\,|x|^{N-|\delta'|}\leq b_N\,\sup_{\{\nu:|\nu|=N-|\delta'|\}}|x^\nu|
\leq  b_N\,\sup_{\{\nu:|\nu|=N\}}|(\partial_x^{\,\delta'} x^\nu)|
\;.
\]
Consequently,
\begin{eqnarray*}
|(\ad\,X)^{\delta'}(Y_j)\,X^{\gamma\backslash\delta'}\varphi|=
|(\partial_x^{\,\delta'}|x|^N)\,\partial_{y_j}\,X^{\gamma\backslash\delta'}
\varphi|
\leq  b_N\,\sup_{\{\nu:|\nu|=N\}}|(\partial_x^{\,\delta'} x^\nu)\,\partial_{y_j}\,X^{\gamma\backslash\delta'}
\varphi|
\end{eqnarray*}
But $x^\nu\partial_{y_j}=X_l$ for a suitable choice of $l\in\{1,\ldots,n+m\}$ and since $\delta'$ has entries in $\{1,\ldots,n\}$
one has
\[
(\partial_x^{\,\delta'} x^\nu)\,\partial_{y_j}\,X^{\gamma\backslash\delta'}
=(\ad\,X)^{\,\delta'}(X_l)\,X^{\gamma\backslash\delta'}
\;.
\]
As  the right hand expression is a sum of monomials of order $|\gamma|+1$ in the $X_j$
one concludes that 
\[
\|(\ad\,X)^{\delta'}(Y_j)\,X^{\gamma\backslash\delta'}\varphi\|_p
\leq 
c_N \sup_{\{\beta:|\beta|=|\gamma|+1\}}\|X^\beta \varphi\|_p\;.
\]
Finally it follows from (\ref{egrrt3.2}) that
\[
\sup_{\{\gamma:|\gamma|=M\}}\|(\ad\, X^\gamma)( Y_j)\varphi\|_p\leq 
c^{\,\prime}_N \sup_{\{\beta:|\beta|=M+1\}}\|X^\beta \varphi\|_p
\]
and this gives the required bound on the second term on the right of (\ref{egrrt3.10}).
This completes the induction.

One concludes that  (\ref{egrrt3.1}) is valid
for all $M\in\{1,\ldots, N+1\}$ and all $p\in\langle1,\infty\rangle$.
Then the boundedness of the $Y$-transforms of all orders up to $N+1$ follows from the boundedness of the $X$-transforms.
The weak-type $(1,1)$ property follows similarly with the aid of Corollary~\ref{cgrrt2.3}.
\hfill$\Box$ 

\bigskip

The $Y$-transforms of order larger than $N+1$ are not necessarily bounded.
Nevertheless the foregoing argument establishes that a large subset of the transforms are bounded.

\begin{cor}\label{cgrrt3.1}
Assume $n\geq 2$.
Let $Y^\alpha=Y^\gamma Y_l$ with $|\alpha|=N+2$.
If either  $l\leq n$ or  $l>n$ and $\max\{\gamma_j:\gamma_j\in \gamma\}>n$   then $Y^\alpha H_N^{-|\alpha|/2}$
extends to a  bounded operator on $L_p(\Ri^{n+m})$
for  each $p\in\langle1,\infty\rangle$ and is of weak-type $(1,1)$.
\end{cor}
\proof\
If $l\leq n$ then $Y_l=X_l$.
Hence  $\|Y^\alpha\varphi\|_p\leq \sup_{\{\gamma:|\gamma|=N+1\}}\|X^\gamma X_l\varphi\|_p$
and the conclusion follows as before

If $l>n$ and  $\max\{\gamma_j:\gamma_j\in \gamma\}>n$ then one uses (\ref{egrrt3.10}) and argues as above. 
The assumption on the indices of $\gamma$ ensures that there are at most $N$ derivatives $\partial_x$ of $|x|^N$ occurring in the estimates. 
Since $n\geq2$ all terms can be bounded as before.
\hfill$\Box$

\bigskip

The $Y$-transforms of order $N+2$ are, however, bounded for small values of~$p$.

\begin{prop}
If $n\geq 2$ the Riesz transforms $Y^\alpha H^{-|\alpha|/2}$ with $|\alpha|= N+2$ extend to bounded operators on 
$L_p(\Ri^{n+m})$ for all $p\in\langle1,n\rangle$ and are weak-type $(1,1)$.

Conversely, if $n\geq1$ then the transforms corresponding to $Y^\alpha=Y^\gamma Y_l$ with   
$\gamma\in\{1,\ldots,n\}^{N+1}$ 
and $l>n$ do not extend to bounded operators on $L_p(\Ri^{n+m})$ for $p\in[n,\infty\rangle$.
\end{prop}
\proof\ 
Consider the first statement.
As a consequence of  Corollary~\ref{cgrrt3.1} it suffices to prove boundedness 
 for $Y^\alpha=Y^\gamma Y_l$ with  $\gamma\in\{1,\ldots,n\}^{N+1}$ and  $l>n$.
But then  $Y^\alpha=X^\gamma Y_l$ and 
\[
Y^\alpha\varphi=Y_l X^\gamma\varphi+{\sum_{\delta\subseteq \gamma}}'
(\ad\, X)^{\delta}( Y_l)X^{\gamma\backslash\delta}\varphi
\]
\vspace{-4mm}

\noindent 
where the prime indicates that the sum is restricted to multi-indices $\delta$ with entries in $\{1,\ldots, n\}$.
Now arguing as in the proof of Theorem~\ref{tgrrt3.2}  the $L_p$-norms of all terms on the right can be bounded
by a multiple of $\sup_{\{\beta:|\beta|=N+2\}}\|X^\beta\varphi\|_p$ with the exception of the leading term in the commutator 
sum, the term $(\ad\, X)^{\gamma}( Y_l)\varphi$.
These latter estimates are valid for all $p\in\langle1,\infty\rangle$.
Now consider the $L_p$-norm of the exceptional term.

Since $|\gamma|=N+1$ and $n\geq 2$ one has
\begin{equation}
\|(\ad\, X)^{\gamma}( Y_l)\varphi\|_p=\|\partial_x^{\gamma}(|x|^N)\partial_{y_l}\varphi\|_p
\leq c_N\,\||x|^{-1}\partial_{y_l}\varphi\|_p
\;.\label{egrrt3.3}
\end{equation}
Next we use the multi-dimensional $L_p$-version of the classical Hardy inequality
(see, for example, \cite{OK} Chapter~2 or \cite{Maz1}, Section~1.3.1).
Explicitly,
\vspace{-2mm}
if $1\leq p<n$ then 
\begin{equation}
\int_{\Ri^n}dx\,|x|^{-p}|\varphi(x)|^p\leq (p/(n-p))^p\int_{\Ri^n}dx\,|(\nabla_x\varphi)(x)|^p
\label{egrrt3.4}
\end{equation}
for all $\varphi\in C^\infty_c(\Ri^n)$.
Combining (\ref{egrrt3.3}) and   (\ref{egrrt3.4})  one finds there is a $c_{p,N}>0$ such that 
\[
\|(\ad\, X)^{\gamma}( Y_l)\varphi\|_p\leq c_{p,N}\,\sup_{1\leq j\leq n}\|\partial_{x_j}\partial_{y_l}\varphi\|_p
\;.
\]
But $\partial_{y_l}=(\ad\, X^\delta)(X_k)$ for a suitable choice $\delta\in\{1,\dots,n\}^N$ and  $k>n$.
In addition $\partial_{x_j}=X_j$.
Therefore one has an estimate
\[
\|(\ad\, X)^{\gamma}( Y_l)\varphi\|_p\leq a_{p,N}\sup_{\{\beta:|\beta|=N+2\}}\|X^\beta\varphi\|_p
\;.
\]
Hence 
\[
\|Y^\alpha\varphi\|_p\leq b_{p,N}\sup_{\{\beta:|\beta|=N+2\}}\|X^\beta\varphi\|_p
\]
with $b_{p,N}>0$.
It then follows from the boundedness of the $X$-transforms on $L_p(\Ri^{n+m})$ with $p\in\langle1,\infty\rangle$ that the $Y^\alpha H^{-|\alpha|/2} $ extend to  bounded operators 
on  the $L_p$-spaces  with $p\in\langle1,n\rangle$.
The weak-type $(1,1)$ property can be established similarly.

Next consider the converse statement.
Since $\gamma\in\{1,\ldots,n\}^{N+1}$ and $Y_j=X_j$ for $j\leq n$ it follows that $Y^\alpha=X^\gamma Y_l$.
But a necessary condition for the transforms to extend to bounded operators is the inclusion 
$C_c^\infty(\Ri^{n+m})\subseteq D(Y^\alpha)$.
Therefore we argue that there is a $\varphi\in C_c^\infty(\Ri^{n+m})$ such that $\varphi\not\in  D(Y^\alpha)$.
Let $\varphi=\varphi_1\varphi_2$ and $\psi=\psi_1\psi_2$ with $\varphi_1,\psi_1\in C_c^\infty(\Ri^n)$ and $\varphi_2,\psi_2\in C_c^\infty(\Ri^m)$.
Then
\[
(-1)^{N+1}(X^\gamma\psi, Y_l\varphi)=(-1)^{N+1}(\partial_x^\gamma\psi_1,|x|^N\varphi_1)\,(\psi_2,\partial_{y_l}\varphi_2)
\;.
\]
Hence for the inclusion  $\varphi\in D(Y^\alpha)$  on $L_p(\Ri^{n+m})$ it is necessary that  $\psi\mapsto(\partial_x^\gamma\psi_1,|x|^N\varphi_1)$ is $L_q$-continuous
where $q$ is the conjugate to $p$.
But $|\gamma|=N+1$ and  if $\varphi>0$ in an open neighbourhood of the origin  $L_q$-continuity requires that $p<n$.
\hfill$\Box$
\bigskip

One  may similarly  deduce  boundedness of higher order $Y$-transforms for a suitable range of $p$ and $n$.
But then the Hardy inequality has to be replaced by the Rellich inequality or an appropriate higher order generalization.
For example, if $|\alpha|=N+3$ then the critical transform to bound is  $Y^\alpha=Y^\gamma Y_l$ with $\gamma\in\{1,\dots,n\}^{N+2}$ and $l>n$.
Then (\ref{egrrt3.3}) is replaced by an estimate
\begin{equation}
\|(\ad\, X)^{\gamma}( Y_l)\varphi\|_p\leq c_N\,\||x|^{-2}\partial_{y_l}\varphi\|_p
\;.\label{egrrt3.31}
\end{equation}
But if $n\geq 3$ and $p\in \langle1,n/2\rangle$ then  the Rellich inequality (see, for example, (3) in  \cite{DaH}) 
states that 
\begin{equation}
\int_{\Ri^n}dx\,|x|^{-2p}|\varphi(x)|^p\leq c_{p,n}^{\,p}\int_{\Ri^n}dx\,|(\Delta_x\varphi)(x)|^p
\label{egrrt3.41}
\end{equation}
for all $\varphi\in C^\infty_c(\Ri^n)$ with $c_{p,n}=p^2/((p-1)n(n-2p))$.
Arguing as above but with (\ref{egrrt3.4}) replaced by (\ref{egrrt3.41}) one concludes the following.

\begin{cor}\label{cgrrt3.2}
If $n\geq 3$ the Riesz transforms $Y^\alpha H^{-|\alpha|/2}$ with $|\alpha|= N+3$ extend to bounded operators on 
$L_p(\Ri^{n+m})$ for all $p\in\langle1,n/2\rangle$ and are weak-type $(1,1)$.
\end{cor}

Again the upper bound $n/2$ is optimal. 
If $p\geq n/2$ there are some transforms of order $N+3$ which are unbounded.

\medskip

The foregoing methods extend to  operators  $-\nabla_{x}^2-c(x)\,\nabla_{y}^2$ with $c$ a sum of positive multiples of even powers of $|x|$. 
But they do not shed light on more general situations such as  $c(x)=|x|^\delta$ with $\delta>0$ but $\delta\neq 2N$.
In particular it would be of interest to understand the properties of Riesz transforms in the case $c(x)=|x|$ analyzed in \cite{ChS}.
It is also not clear if the $|x|^{2N}$-results for  are still valid if one only has $c(x)\sim |x|^{2N}$, i.e. if 
$a\,|x|^{2N}\leq c(x)\leq b\,|x|^{2N}$
for some $a,b>0$, although many properties of Gru\v{s}in operators are known to be invariant under such equivalence relations
\cite{RSi2} \cite{RSi6a}.

\end{document}